\documentclass[12pt,twoside]{article}
\usepackage{amsmath, amssymb, amsthm,graphicx}
\font\fourteenb=cmb10 at 14pt
\begin{document}
\vspace*{-1.0cm}\noindent \copyright
 Journal of Technical University at Plovdiv\\[-0.0mm]\
\ Fundamental Sciences and Applications, Vol. 2, 1996\\[-0.0mm]
\textit{Series A-Pure and Applied Mathematics}\\[-0.0mm]
\ Bulgaria, ISSN 1310-8271\\[+1.2cm]
\font\fourteenb=cmb10 at 12pt
\begin{center}

{\bf \LARGE On some interpolation theorems for the multipliers of
the Cauchy- Stiltjes type integrals
   \\ \ \\ \large Peyo Stoilov}
\end{center}

\footnotetext{{\bf 1991 Mathematics Subject Classification:}
Primary 30E20, 30D50.} \footnotetext{{\it Key words and phrases:}
Free interpolation, Cauchy - Stiltjes type integral, multipliers.
}

\begin{abstract}
Let  $m_{\alpha } {\kern 1pt} $ ($\alpha >0$)  denote the set of
all multipliers of the analytic functions in the unit disk D,
representable  by a Cauchy - Stiltjes type integral  $\int \limits
_{T} 1/(1-\overline{\zeta }z)^{\alpha } {\kern 1pt} d\mu (\zeta
).$

Let  the  sequence   $a=\{ a_{k} \} _{k\ge 1} \subset D{\kern 1pt}
$  and

$$R_{a} f=\{ f(a_{k} )\} _{k\ge 1} {\kern 1pt} ,{\kern 1pt} {\kern 1pt} {\kern 1pt} {\kern 1pt} {\kern 1pt} {\kern 1pt} {\kern 1pt} {\kern 1pt} {\kern 1pt} {\kern 1pt} {\kern 1pt} {\kern 1pt} {\kern 1pt} {\kern 1pt} {\kern 1pt} {\kern 1pt} {\kern 1pt} {\kern 1pt} {\kern 1pt} {\kern 1pt} {\kern 1pt} {\kern 1pt} {\kern 1pt} {\kern 1pt} {\kern 1pt} {\kern 1pt} {\kern 1pt} {\kern 1pt} {\kern 1pt} S_{a} f=\{ f'(a_{k} )(1-\left|a_{k} \right|^{2} \} _{k\ge 1} .{\kern 1pt} $$

In this paper we prove that  if  a  sequence   $a=\{ a_{k} \}
_{k\ge 1} \subset D{\kern 1pt} $  satisfies conditions of  Newman
- Carleson  and  Stoltz {\it    }, then

$$R_{a} f=b_{v} {\kern 1pt} ,{\kern 1pt} {\kern 1pt} {\kern 1pt} {\kern 1pt} {\kern 1pt} {\kern 1pt} S_{a} f=l^{1} .{\kern 1pt} $$
\end{abstract}

     \section{Introduction}

     Let  $D{\kern 1pt} $  denote the unit disk in the complex plane and   $T{\kern 1pt} $ - the unit circle. For  $0<p\le \infty $   let   $H^{P} $   be the usual Hardy class [1].

     Let  $M{\kern 1pt} $  be the Banach space of all complex-valued Borel measures on  $T{\kern 1pt} $  with the usual variation norm. For  $\alpha >0$ ,  let  $F_{\alpha } $  denote the family of all functions  $g{\kern 1pt} $  for which there exists  $\mu \in M{\kern 1pt} $  such that
\begin{equation}
    g(z)=\int \limits _{T} \frac{1}{(1-\overline{\zeta }z)^{\alpha } } {\kern 1pt} d\mu (\zeta ),\quad z\in
     D.
\end{equation}

     We note that  $F_{\alpha } {\kern 1pt} $  is a Banach space with the natural norm

     $${\rm \parallel }g{\rm \parallel }_{F_{\alpha } } =\inf \left\{{\rm \parallel }\mu {\rm \parallel }:{\kern 1pt} {\kern 1pt} {\kern 1pt} {\kern 1pt} {\kern 1pt} {\kern 1pt} \mu \in M{\kern 1pt} {\kern 1pt} {\kern 1pt} {\kern 1pt} such{\kern 1pt} {\kern 1pt} {\kern 1pt} that{\kern 1pt} {\kern 1pt} {\kern 1pt} {\kern 1pt} (1){\kern 1pt} {\kern 1pt} {\kern 1pt} {\kern 1pt} {\kern 1pt} holds\; \right\}.$$

     The family  $F_{\alpha } {\kern 1pt} $  was introduced in [2] and the following results were proved
     there:

     If  $f\in F_{\alpha } {\kern 1pt} $ ,  $g\in F_{\beta } {\kern 1pt} $ ,  then   $f{\kern 1pt} g\in F_{\alpha +\beta } {\kern 1pt} $    $(\alpha >0,{\kern 1pt} {\kern 1pt} {\kern 1pt} {\kern 1pt} \beta >0{\kern 1pt} {\kern 1pt} )$ .

\

     If  $\alpha <\beta ,{\kern 1pt} $  then  $F_{\alpha } \subset F_{\beta } {\kern 1pt} $ .

\

     {\bf Definition.}  {\it Suppose that } $f{\kern 1pt} $ {\it is holomorphic in } $D$ {\it . Then } $f{\kern 1pt} $ {\it  is called }{\it  }{\it a }

      {\it multiplier of } $F_{\alpha } {\kern 1pt} $ {\it  if }{\it  } $g\in F_{\alpha } \Rightarrow fg\in F_{\alpha } $ {\it .}

\

     Let  $m_{\alpha } {\kern 1pt} $  denote the set of  all multipliers of  $F_{\alpha } {\kern 1pt} $  and
     $${\rm \parallel }f{\rm \parallel }_{m_{\alpha } } =\sup \left\{{\rm \parallel }fg{\rm \parallel }_{F_{\alpha } } :{\rm \parallel }g{\rm \parallel }_{F_{\alpha } } \le 1\right\}.$$

      Various properties of  $m_{\alpha } {\kern 1pt} $ were studied in [3], [8]. For example  it is proved that  $\alpha <\beta $  implies   $m_{\alpha } {\kern 1pt} \subset m_{\alpha } \subset H^{\infty } $ .  Also it is proved that  if   $f'\in H^{1} $,   then  $f\in m_{\alpha } $  for all  $\alpha >0$ .

     These results will be applied here.

     In this paper some interpolation theorems for  $m_{1} {\kern 1pt} $  and  $F_{1} {\kern 1pt} $  obtained in [4] and [7] are generalized for  $m_{\alpha } {\kern 1pt} $ and   $F_{\alpha } {\kern 1pt} $ .

     \section{Interpolation theorems}

\

     {\bf Definition.}{\it We say that a sequence } $a=\{ a_{k}
\} _{k\ge 1} \subset D{\kern 1pt} $ {\it  satisfies Newman - }{\it
Carleson condition (condition (} $N-C$ {\it )) if }

$$\delta (a)=\inf {\kern 1pt} {\kern 1pt} {\kern 1pt} ({\kern 1pt} {\kern 1pt} {\kern 1pt} {\kern 1pt} {\kern 1pt} {\kern 1pt} {\kern 1pt} \prod \limits _{k\ne n} \left|\frac{a_{k} -a_{n} }{1-\overline{a}_{k} a_{n} } \right|:\; n=1,2,\ldots ){\kern 1pt} {\kern 1pt} {\kern 1pt} {\kern 1pt} {\kern 1pt} {\kern 1pt} {\kern 1pt} >0.$$

     {\bf Definition.}{\it We say that a sequence } $a=\{ a_{k} \} _{k\ge 1} \subset D{\kern 1pt} $ {\it  satisfies the Stoltz condition (condition (} $S$ {\it ))  if }

$$\Lambda (a)=\sup {\kern 1pt} {\kern 1pt} {\kern 1pt} {\kern 1pt} ({\kern 1pt} {\kern 1pt} {\kern 1pt} {\kern 1pt} {\kern 1pt} {\kern 1pt} \frac{\left|1-a_{k} \right|}{1-\left|a_{k} \right|} :{\kern 1pt} {\kern 1pt} {\kern 1pt} \; k=1,2,\ldots ){\kern 1pt} {\kern 1pt} {\kern 1pt} {\kern 1pt} {\kern 1pt} <\infty .$$

     We will need the following Lemma.

     {\bf Lemma 1.} {\it If a sequence } $a=\{ a_{k} \} _{k\ge 1} \subset D{\kern 1pt} $ {\it  satisfy conditions (} $N-C$ {\it ) and (} $S$ {\it ) then}
$$\Omega (a){\kern 1pt} {\kern 1pt} {\kern 1pt} {\kern 1pt} ={\kern 1pt} {\kern 1pt} {\kern 1pt} \mathop{\sup }\limits_{y\ge 0} {\kern 1pt} {\kern 1pt} {\kern 1pt} {\kern 1pt} \sum \limits _{k\ge 1} \frac{(1-|a_{k} |)y}{(1-|a_{k} |)^{2} +y^{2} } <\infty .$$

     The proof of this Lemma can be found, for example  in [5] (10.2 Lemma 11) or [6].

     Let  $l^{1} $  and  $b_{v} $ is the class of sequences:
     $$\displaystyle l^{1} ={\kern 1pt} {\kern 1pt} {\kern 1pt} ({\kern 1pt} x=\{ x_{k} \} _{k\ge 1} :{\kern 1pt} {\kern 1pt} {\kern 1pt} {\kern 1pt} {\kern 1pt} {\kern 1pt} {\kern 1pt} {\kern 1pt} \left\| x\right\| _{l^{1} } =\sum \limits _{k\ge 1}\left|x_{k} \right| {\kern 1pt} {\kern 1pt} {\kern 1pt} <{\kern 1pt} {\kern 1pt} {\kern 1pt} \infty {\kern 1pt} {\kern 1pt} {\kern 1pt} {\kern 1pt} ),$$

     $$b_{v} ={\kern 1pt} {\kern 1pt} {\kern 1pt} ({\kern 1pt} x=\{ x_{k} \} _{k\ge 1} :{\kern 1pt} {\kern 1pt} {\kern 1pt} {\kern 1pt} {\kern 1pt} {\kern 1pt} {\kern 1pt} {\kern 1pt} \left\| x\right\| _{b_{v} } =\sum \limits _{k\ge 1}\left|x_{k+1} -x_{k} \right| {\kern 1pt} {\kern 1pt} {\kern 1pt} <{\kern 1pt} {\kern 1pt} {\kern 1pt} \infty {\kern 1pt} {\kern 1pt} ).$$

     The following two theorems were proved in [4] for  $\alpha =1.$

     \

{\bf Theorem 1.}  {\it Let } $\alpha >0$ {\it  and the sequence }
$a=\{ a_{k} \} _{k\ge 1} {\kern 1pt} $ {\it  } $({\kern 1pt}
{\kern 1pt} a_{k} \ne a_{n} {\kern 1pt} {\kern 1pt} {\kern 1pt}
{\kern 1pt} {\kern 1pt} if{\kern 1pt} {\kern 1pt} {\kern 1pt}
{\kern 1pt} {\kern 1pt} k\ne n{\kern 1pt} {\kern 1pt} {\kern 1pt}
{\kern 1pt} {\kern 1pt} {\kern 1pt} {\kern 1pt} and{\kern 1pt}
{\kern 1pt} {\kern 1pt} {\kern 1pt} {\kern 1pt} {\kern 1pt}
\left|a_{k} \right|\le \left|a_{k+1} \right|{\kern 1pt} {\kern
1pt} {\kern 1pt} )$ {\it  }{\it satisfy conditions (} $N-C$ {\it )
and (} $S$ {\it ). }

{\it Then}

     {\it a) }{\it If  } $f\in m_{\alpha } ,$  {\it then  } $\{ f(a_{k} )\} _{k\ge 1} {\kern 1pt} \in b_{v} ;{\kern 1pt} $

     {\it b)} {\it For each  sequence }{\it  } $x=\{ x_{k} \} _{k\ge 1} \in b_{v} $ {\it , there is a function } $f$ {\it  in } $m_{\alpha } $ {\it  ,
     such that }{\it   }    $f(a_{k} )=x_{k} ,{\kern 1pt} {\kern 1pt} {\kern 1pt} {\kern 1pt} {\kern 1pt} {\kern 1pt} {\kern 1pt} {\kern 1pt} k=1,2,\ldots .$

\

     {\it Proof.} {\it a)} Let  $f\in m_{\alpha } .$  Then

\

      ${\kern 1pt} \left\| \{ f(a_{k} )\} \right\| _{b_{v} } =\sum \limits _{k\ge 1}\left|f(a_{k} )-f(a_{k+1} )\right| {\kern 1pt} {\kern 1pt} {\kern 1pt} \le $  $\sum \limits _{k\ge 1}\left|f(a_{k} )-f(\left|a_{k} \right|)\right|
      +$

      \

      $+\sum \limits _{k\ge 1}\left|f(\left|a_{k} \right|)-f(\left|a_{k+1} \right|)\right| +\sum \limits _{k\ge 1}\left|f(\left|a_{k+1} \right|-f(a_{k+1} )\right| \le $  $2\Lambda (f){\kern 1pt} {\kern 1pt} {\kern 1pt} {\kern 1pt} {\kern 1pt} {\kern 1pt} +{\kern 1pt} {\kern 1pt} {\kern 1pt} {\kern 1pt} I(f)$ {\it ,}

\

     where

              $$I(f)=\sum \limits _{k\ge 1}\left|f(\left|a_{k} \right|)-f(\left|a_{k+1} \right|)\right| ,$$

              $$\Lambda (f){\kern 1pt} {\kern 1pt} {\kern 1pt} {\kern 1pt} {\kern 1pt} {\kern 1pt} =\sum \limits _{k\ge 1}\left|f(a_{k} )-f(\left|a_{k} \right|)\right| .$$

     In [3] (Theorem 2.6) it was proved that   $\int _{0}^{1}\left|f'(r)\right| dr<\infty $   and consequently

$$I(f)\le \sum \limits _{k\ge 1}\int \limits _{\left|a_{k+1} \right|}^{\left|a_{k} \right|}\left|f'(r)\right|dr  <\int _{0}^{1}\left|f'(r)\right| dr<\infty .$$

     We will show that  $\Lambda (f)<\infty $  too.

     Let  $n\ge [\alpha ]+1$   is a natural number. Then  $m_{\alpha } {\kern 1pt} \subset m_{n} $  and  $f\in m_{n} .$

     Since   $\left\| 1/(1-z)^{n} \right\| _{F_{n} {\kern 1pt} } =1$ , then there exists a measurer   $\mu \in M$  for which

     $$\frac{f(z)}{(1-z)^{n} } =\int \limits _{T} \frac{1}{(\zeta -z)^{n} } {\kern 1pt} d\mu (\zeta )  {\kern 1pt} {\kern 1pt}   \Rightarrow  {\kern 1pt} {\kern 1pt}    f(z)=\int \limits _{T} \left(\frac{1-z}{\zeta -z} \right)^{n} {\kern 1pt} d\mu (\zeta
                     ).$$

     For  $\Lambda (f)$  we obtain

     $$\Lambda (f){\kern 1pt} {\kern 1pt} {\kern 1pt} {\kern 1pt} {\kern 1pt} {\kern 1pt} =\sum \limits _{k\ge 1}\left|\int \limits _{T}\left(\frac{1-a_{k} }{\zeta -a_{k} } \right)^{n} -\left(\frac{1-\left|a_{k} \right|}{\zeta -\left|a_{k} \right|} \right)^{n} {\kern 1pt} {\kern 1pt} {\kern 1pt} {\kern 1pt} d\mu (\zeta ) \right| {\kern 1pt} {\kern 1pt} {\kern 1pt} {\kern 1pt} {\kern 1pt} \le $$
     $$\le \sum \limits _{k\ge 1}\int \limits _{T}\left|\frac{1-a_{k} }{\zeta -a_{k} } -\frac{1-\left|a_{k} \right|}{\zeta -\left|a_{k} \right|} \right|{\kern 1pt} {\kern 1pt} {\kern 1pt} {\kern 1pt} L(\zeta ){\kern 1pt} {\kern 1pt} {\kern 1pt} d\left|\mu \right|(\zeta )  {\kern 1pt} {\kern 1pt} {\kern 1pt} {\kern 1pt} {\kern 1pt} ={\kern 1pt} $$
     $$=\sum \limits _{k\ge 1}\int \limits _{T}\frac{\left|1-\varsigma \right|\left|a_{k} -\left|a_{k} \right|\right|}{\left|\zeta -a_{k} \right|\left|\zeta -\left|a_{k} \right|\right|} {\kern 1pt} {\kern 1pt} L(\zeta ){\kern 1pt} {\kern 1pt} {\kern 1pt} d\left|\mu \right|(\zeta ),$$

where

     $$L(\zeta )=\sum \limits _{k=0}^{n-1}\left|\frac{1-a_{k} }{\zeta -a_{k} } \right| ^{k} \left|\frac{1-\left|a_{k} \right|}{\zeta -\left|a_{k} \right|} \right|^{n-k-1} {\kern 1pt} {\kern 1pt} {\kern 1pt} {\kern 1pt} {\kern 1pt} \le $$

     $$\le {\kern 1pt} {\kern 1pt} {\kern 1pt} {\kern 1pt} {\kern 1pt} \sum \limits _{k=0}^{n-1}\left|\frac{1-a_{k} }{1-\left|a_{k} \right|} \right| ^{k} {\kern 1pt} {\kern 1pt} {\kern 1pt} {\kern 1pt} \le \sum \limits _{k=0}^{n-1}\lambda ^{k}  {\kern 1pt} {\kern 1pt} {\kern 1pt} ={\kern 1pt} {\kern 1pt} \frac{1-\lambda ^{n} }{1-\lambda } .$$

     Therefore

     $$\Lambda (f){\kern 1pt} {\kern 1pt} {\kern 1pt} {\kern 1pt} {\kern 1pt} {\kern 1pt} \le {\kern 1pt} {\kern 1pt} \frac{1-\lambda ^{n} }{1-\lambda } {\kern 1pt} {\kern 1pt} \sum \limits _{k\ge 1}\int \limits _{T}\frac{\left|1-\zeta \right|\left|a_{k} -\left|a_{k} \right|\right|}{\left|\zeta -a_{k} \right|\left|\zeta -\left|a_{k} \right|\right|} {\kern 1pt} {\kern 1pt} {\kern 1pt} {\kern 1pt} {\kern 1pt} d\left|\mu \right|(\zeta )  {\kern 1pt} .$$

     Using the condition (S), we obtain
\begin{equation}
\frac{\left|a_{k} -\left|a_{k} \right|\right|}{1-\left|a_{k}
\right|} =\frac{\left|a_{k} -1+1-\left|a_{k}
\right|\right|}{1-\left|a_{k} \right|} {\kern 1pt} {\kern 1pt}
{\kern 1pt} {\kern 1pt} \le {\kern 1pt} {\kern 1pt} {\kern 1pt}
{\kern 1pt} {\kern 1pt} 1+\lambda ,
\end{equation}
\begin{equation}
\left|\frac{\zeta -\left|a_{k} \right|}{\zeta -a_{k} } \right|\le
1+\frac{\left|a_{k} -\left|a_{k} \right|\right|}{\left|\zeta
-a_{k} \right|} {\kern 1pt} {\kern 1pt} {\kern 1pt} {\kern 1pt}
\le {\kern 1pt} {\kern 1pt} {\kern 1pt} 1+\frac{\left|a_{k}
-\left|a_{k} \right|\right|}{1-\left|a_{k} \right|} {\kern 1pt}
{\kern 1pt} \le {\kern 1pt} {\kern 1pt} {\kern 1pt} {\kern 1pt}
2+\lambda .
\end{equation}

\

     It follows from the above inequalities that

     $$\Lambda (f){\kern 1pt} {\kern 1pt} {\kern 1pt} {\kern 1pt} {\kern 1pt} {\kern 1pt} \le {\kern 1pt} {\kern 1pt} \frac{1-\lambda ^{n} }{1-\lambda } {\kern 1pt} {\kern 1pt} (1+\lambda )(2+\lambda )\int \limits _{T}{\kern 1pt} {\kern 1pt} \sum \limits _{k\ge 1}\frac{\left|1-\zeta \right|({\kern 1pt} {\kern 1pt} 1-\left|a_{k} \right|{\kern 1pt} {\kern 1pt} )}{\left|\zeta -\left|a_{k} \right|\right|^{2} }   {\kern 1pt} {\kern 1pt} {\kern 1pt} {\kern 1pt} {\kern 1pt} {\kern 1pt} d\left|\mu \right|(\zeta ){\kern 1pt} {\kern 1pt} {\kern 1pt} \le $$

     $${\kern 1pt} {\kern 1pt} {\kern 1pt} {\kern 1pt} {\kern 1pt} \le {\kern 1pt} {\kern 1pt} \frac{1-\lambda ^{n} }{1-\lambda } {\kern 1pt} {\kern 1pt} (1+\lambda )(2+\lambda )\left\| \mu \right\| {\kern 1pt} {\kern 1pt} {\kern 1pt} {\kern 1pt} {\kern 1pt} {\kern 1pt} \mathop{\sup }\limits_{\zeta \in T} {\kern 1pt} {\kern 1pt} {\kern 1pt} {\kern 1pt} \sum \limits _{k\ge 1}\frac{\left|1-\zeta \right|({\kern 1pt} {\kern 1pt} 1-\left|a_{k} \right|{\kern 1pt} {\kern 1pt} )}{\left|\zeta -\left|a_{k} \right|\right|^{2} }  \le $$

     $${\kern 1pt} {\kern 1pt} {\kern 1pt} {\kern 1pt} {\kern 1pt} \le \frac{1-\lambda ^{n} }{1-\lambda } {\kern 1pt} {\kern 1pt} (1+\lambda )(2+\lambda )\left\| \mu \right\| {\kern 1pt} {\kern 1pt} {\kern 1pt} {\kern 1pt} {\kern 1pt} {\kern 1pt} \mathop{\sup }\limits_{\zeta \in T} {\kern 1pt} {\kern 1pt} {\kern 1pt} {\kern 1pt} \sum \limits _{k\ge 1}\frac{\left|1-\zeta \right|({\kern 1pt} {\kern 1pt} 1-\left|a_{k} \right|{\kern 1pt} {\kern 1pt} )}{(1-|a_{k} |)^{2} +|a_{1} |^{2} \left|1-\zeta \right|^{2} } . $$

     Here we used the inequalities

     $$\left|\zeta -\left|a_{k} \right|\right|^{2} =(\zeta -\left|a_{k} \right|)(\overline{\zeta }-\left|a_{k} \right|)=1-\left|a_{k} \right|(\zeta +\overline{\zeta })+\left|a_{k} \right|^{2} =(1-\left|a_{k} \right|)^{2} +\left|a_{k} \right|({\kern 1pt} {\kern 1pt} 2-\zeta -\overline{\zeta }{\kern 1pt} {\kern 1pt} )=$$

     $$=(1-\left|a_{k} \right|)^{2} +\left|a_{k} \right|(1-\zeta )(1-\overline{\zeta })=(1-\left|a_{k} \right|)^{2} +\left|a_{k} \right|\left|1-\zeta \right|^{2} \ge (1-\left|a_{k} \right|)^{2} +\left|a_{k} \right|^{2} \left|1-\zeta \right|^{2} \ge $$

     $$\ge (1-\left|a_{k} \right|)^{2} +\left|a_{1} \right|^{2} \left|1-\zeta \right|^{2} .$$

     Applying Lemma 1 for  $y=\left|a_{1} \right|\left|1-\zeta \right|$ ,  we  obtain

$${\kern 1pt} {\kern 1pt} {\kern 1pt} {\kern 1pt} \Lambda (f){\kern 1pt} {\kern 1pt} {\kern 1pt} {\kern 1pt} {\kern 1pt} {\kern 1pt} \le {\kern 1pt} {\kern 1pt} {\kern 1pt} \left\| f\right\| _{m_{n} } \frac{1-\lambda ^{n} }{1-\lambda } {\kern 1pt} {\kern 1pt} (1+\lambda )(2+\lambda ){\kern 1pt} {\kern 1pt} {\kern 1pt} \Omega (\{ f(a_{k} )\} )/\left|a_{1} \right|{\kern 1pt} {\kern 1pt} {\kern 1pt} {\kern 1pt} {\kern 1pt} {\kern 1pt} {\kern 1pt} <\infty .$$

     Consequently

\

      ${\kern 1pt} \left\| \{ f(a_{k} )\} \right\| _{b_{v} } {\kern 1pt} {\kern 1pt} {\kern 1pt} \le 2\Lambda (f){\kern 1pt} {\kern 1pt} {\kern 1pt} {\kern 1pt} {\kern 1pt} {\kern 1pt} +{\kern 1pt} {\kern 1pt} {\kern 1pt} {\kern 1pt} I(f)<\infty $ {\it  } and{\it  } $\{ f(a_{k} )\} _{k\ge 1} \in b_{v} $ .

\

     {\it b)} Let  $x=\{ x_{k} \} _{k\ge 1} \in b_{v} .$ {\it  }From the theorem of  Vinogradov [6] it follows that there exists an analytic function  $f$ ,  such that{\it  } $f'\in H^{1} $ {\it  }and  $f(a_{k} )=x_{k} ,{\kern 1pt} {\kern 1pt} {\kern 1pt} {\kern 1pt} {\kern 1pt} {\kern 1pt} {\kern 1pt} {\kern 1pt} k=1,2,\ldots .$  Since  $f'\in H^{1} $ {\it ,  }it follows that  $f\in m_{\alpha } $   for all  $\alpha >0$  [3].

\

{\bf Theorem 2.} {\it Let }$\alpha >0$ {\it  and the sequence }
$a=\{ a_{k} \} _{k\ge 1} {\kern 1pt} $ {\it  } $({\kern 1pt}
{\kern 1pt} a_{k} \ne a_{n} {\kern 1pt} {\kern 1pt} {\kern 1pt}
{\kern 1pt} {\kern 1pt} if{\kern 1pt} {\kern 1pt} {\kern 1pt}
{\kern 1pt} {\kern 1pt} k\ne n{\kern 1pt} {\kern 1pt} {\kern 1pt}
{\kern 1pt} {\kern 1pt} {\kern 1pt} {\kern 1pt} and{\kern 1pt}
{\kern 1pt} {\kern 1pt} {\kern 1pt} {\kern 1pt} {\kern 1pt}
\left|a_{k} \right|\le \left|a_{k+1} \right|{\kern 1pt} {\kern
1pt} {\kern 1pt} )$ {\it satisfy conditions (} $N-C$ {\it ) and (}
$S$ {\it ). }

{\it Then}

     {\it a) }{\it If  } $g\in F_{\alpha } ,$  {\it then  } $\{ (1-a_{k} )^{\alpha } {\kern 1pt} {\kern 1pt} g(a_{k} )\} _{k\ge 1} {\kern 1pt} {\kern 1pt} \in b_{v} ;{\kern 1pt} $

     {\it b)} {\it For each }{\it  }{\it sequence } $x=\{ x_{k} \} _{k\ge 1} \in b_{v} $ {\it , there is a function } $g$ {\it  in } $F_{\alpha } $ {\it ,  such that }{\it  } $(1-a_{k} )^{\alpha } g(a_{k} )=x_{k} ,{\kern 1pt} {\kern 1pt} {\kern 1pt} {\kern 1pt} {\kern 1pt} {\kern 1pt} {\kern 1pt} {\kern 1pt} k=1,2,\ldots .$

\

     {\it Proof.} {\it a)} Let  $g\in F_{\alpha } $  and $f(z)=(1-z)^{\alpha } {\kern 1pt} {\kern 1pt} g(z)$ . If  $n\ge [\alpha ]+1$   is a natural number, then

     $$\frac{f(z)}{(1-z)^{n} } =\frac{f(z)}{(\zeta -z)^{\alpha } } \cdot \frac{1}{(\zeta -z)^{n-\alpha } } =g(z)\cdot \frac{1}{(\zeta -z)^{n-\alpha } } \in F_{\alpha +n-\alpha } =F_{n} $$

and  from the proof of  Theorem 1 follows

     $$\{ f(a_{k} )\} _{k\ge 1} =\{ (1-a_{k} )^{\alpha } {\kern 1pt} {\kern 1pt} g(a_{k} )\} _{k\ge 1} {\kern 1pt} {\kern 1pt} \in b_{v} .$$

\

     {\it b)} Let  $x=\{ x_{k} \} _{k\ge 1} \in b_{v} .$  From Theorem 1 it follows that there exists a function   $f\in m_{\alpha } $ {\it  }such that{\it  } $f(a_{k} )=x_{k} ,{\kern 1pt} {\kern 1pt} {\kern 1pt} {\kern 1pt} {\kern 1pt} {\kern 1pt} {\kern 1pt} {\kern 1pt} k=1,2,....$ Then  $g(z)=\frac{f(z)}{(1-z)^{\alpha } } \in F_{\alpha } $  and   $(1-a_{k} )^{\alpha } g(a_{k} )=x_{k} ,{\kern 1pt} {\kern 1pt} {\kern 1pt} {\kern 1pt} {\kern 1pt} {\kern 1pt} {\kern 1pt} {\kern 1pt} k=1,2,\ldots .$

\

     The following theorem generalizes a result in [7].

\

{\bf Theorem 3.} {\it Let } $\alpha >0$ {\it  and the sequence }
$a=\{ a_{k} \} _{k\ge 1} {\kern 1pt} $ {\it  satisfy conditions (}
$N-C$ {\it ) and (} $S$ {\it ). }

{\it Then}

     {\it a) }{\it If  } $f\in m_{\alpha } ,$  {\it then  } $\{ {\kern 1pt} {\kern 1pt} {\kern 1pt} f'(a_{k} )(1-\left|a_{k} \right|^{2} {\kern 1pt} {\kern 1pt} \} _{k\ge 1} {\kern 1pt} \in l_{1} ;{\kern 1pt} $

     {\it b)} {\it For each  sequence } $x=\{ x_{k} \} _{k\ge 1} \in l_{1} $ {\it , there is a function } $f$ {\it  in } $m_{\alpha } $ {\it   }{\it such that }{\it  } $f'(a_{k} )(1-\left|a_{k} \right|^{2} )=x_{k} ,{\kern 1pt} {\kern 1pt} {\kern 1pt} {\kern 1pt} {\kern 1pt} {\kern 1pt} {\kern 1pt} {\kern 1pt} k=1,2,\ldots .$

\

     {\it Proof.} {\it a)} Let  $f\in m_{\alpha } $  and  $n\ge [\alpha ]+1$   is a natural number. Since $f\in m_{n} $  then there exists a measure  $\mu \in M$  for which

     $$\frac{f(z)}{(1-z)^{n} } =\int \limits _{T} \frac{1}{(\zeta -z)^{n} } {\kern 1pt} d\mu (\zeta );$$

     $$f(z)=\int \limits _{T} \left(\frac{1-z}{\zeta -z} \right)^{n} {\kern 1pt} d\mu (\zeta );$$

     $$f'(z)=n\int \limits _{T} \left(\frac{1-z}{\zeta -z} \right)^{n-1} {\kern 1pt} \frac{1-\zeta }{(\zeta -z)^{2} } d\mu (\zeta ).$$

     Using the condition (S)  and (3) ,  we obtain

     $${\kern 1pt} {\kern 1pt} \sum \limits _{k\ge 1}\left|f'(a_{k} )\right|(1-\left|a_{k} \right|^{2} )\le n\int \limits _{T}\sum \limits _{k\ge 1}{\kern 1pt} {\kern 1pt} \left|\frac{1-a_{k} }{\zeta -a_{k} } \right|^{n-1} {\kern 1pt} {\kern 1pt} {\kern 1pt} {\kern 1pt} \frac{\left|1-\zeta \right|(1-\left|a_{k} \right|^{2} )}{\left|\zeta -a_{k} \right|^{2} } {\kern 1pt} {\kern 1pt} {\kern 1pt} {\kern 1pt} {\kern 1pt} {\kern 1pt} d\left|\mu \right|(\zeta )   {\kern 1pt} {\kern 1pt} {\kern 1pt} {\kern 1pt} {\kern 1pt} \le $$

     $${\kern 1pt} {\kern 1pt} {\kern 1pt} {\kern 1pt} {\kern 1pt} \le 2n\lambda ^{n-1} (2+\lambda )^{2} \int \limits _{T}\sum \limits _{k\ge 1}{\kern 1pt} {\kern 1pt} {\kern 1pt} {\kern 1pt} {\kern 1pt} \frac{\left|1-\zeta \right|(1-\left|a_{k} \right|)}{\left|\zeta -\left|a_{k} \right|\right|^{2} } {\kern 1pt} {\kern 1pt} {\kern 1pt} {\kern 1pt} {\kern 1pt} {\kern 1pt} d\left|\mu \right|(\zeta )  $$

     $${\kern 1pt} {\kern 1pt} {\kern 1pt} {\kern 1pt} {\kern 1pt} \le 2n\lambda ^{n-1} (2+\lambda )^{2} \left\| \mu \right\| {\kern 1pt} {\kern 1pt} {\kern 1pt} {\kern 1pt} {\kern 1pt} {\kern 1pt} \mathop{\sup }\limits_{\zeta \in T} {\kern 1pt} {\kern 1pt} {\kern 1pt} {\kern 1pt} \sum \limits _{k\ge 1}\frac{\left|1-\zeta \right|({\kern 1pt} {\kern 1pt} 1-\left|a_{k} \right|{\kern 1pt} {\kern 1pt} )}{\left|\zeta -\left|a_{k} \right|\right|^{2} }  $$

     $${\kern 1pt} {\kern 1pt} {\kern 1pt} {\kern 1pt} {\kern 1pt} \le 2n\lambda ^{n-1} (2+\lambda )^{2} \left\| \mu \right\| {\kern 1pt} {\kern 1pt} {\kern 1pt} {\kern 1pt} {\kern 1pt} {\kern 1pt} \mathop{\sup }\limits_{\zeta \in T} {\kern 1pt} {\kern 1pt} {\kern 1pt} {\kern 1pt} \sum \limits _{k\ge 1}\frac{\left|1-\zeta \right|({\kern 1pt} {\kern 1pt} 1-\left|a_{k} \right|{\kern 1pt} {\kern 1pt} )}{(1-|a_{k} |)^{2} +|a_{1} |^{2} \left|1-\zeta \right|^{2} } . $$

     Applying Lemma 1 for  $y=\left|a_{1} \right|\left|1-\zeta \right|$ ,  we   obtain

$${\kern 1pt} {\kern 1pt} {\kern 1pt} {\kern 1pt} \sum \limits _{k\ge 1}\left|f'(a_{k} )\right|(1-\left|a_{k} \right|^{2} ) {\kern 1pt} {\kern 1pt} {\kern 1pt} {\kern 1pt} {\kern 1pt} {\kern 1pt} \le {\kern 1pt} {\kern 1pt} {\kern 1pt} \left\| f\right\| _{m_{n} } 2n\lambda ^{n-1} (2+\lambda )^{2} {\kern 1pt} {\kern 1pt} {\kern 1pt} \Omega (\{ f(a_{k} )\} )/\left|a_{1} \right|{\kern 1pt} {\kern 1pt} {\kern 1pt} {\kern 1pt} {\kern 1pt} {\kern 1pt} {\kern 1pt} <\infty .$$

     Consequently  $\{ f(a_{k} )\} _{k\ge 1} \in l_{1} .$

\

     {\it b)} Let  $x=\{ x_{k} \} _{k\ge 1} \in l_{1} .$  We   construct the following function

     $$h(z)=\sum \limits _{n\ge 1}\frac{(1-\left|a_{n} \right|^{2} )B_{n} (z)}{(1-{\kern 1pt} {\kern 1pt} \overline{a}_{n} {\kern 1pt} z)^{2} B_{n} (a_{n} )}  {\kern 1pt} x_{n} {\kern 1pt} {\kern 1pt} ,{\kern 1pt} {\kern 1pt} {\kern 1pt} {\kern 1pt} {\kern 1pt} {\kern 1pt} {\kern 1pt} {\kern 1pt} {\kern 1pt} {\kern 1pt} {\kern 1pt} {\kern 1pt} {\kern 1pt} {\kern 1pt} {\kern 1pt} where{\kern 1pt} {\kern 1pt} {\kern 1pt} {\kern 1pt} {\kern 1pt} {\kern 1pt} {\kern 1pt} {\kern 1pt} {\kern 1pt} {\kern 1pt} {\kern 1pt} B_{n} (z)=\prod \limits _{\mathop{k\ne n}\limits^{k=1} }^{\infty }\frac{z-a_{k} }{1-{\kern 1pt} {\kern 1pt} \overline{a}_{k} {\kern 1pt} z}  \cdot \frac{\left|a_{k} \right|}{a_{k} } {\kern 1pt} {\kern 1pt} .$$

     Since

     $$\left\| \frac{1-\left|a_{n} \right|^{2} }{(1-{\kern 1pt} {\kern 1pt} \overline{a}_{n} {\kern 1pt} z)^{2} } \right\| _{H^{1} } =\frac{1}{2\pi } \int \limits _{T}\frac{1-\left|a_{n} \right|^{2} }{\left|1-{\kern 1pt} {\kern 1pt} \overline{a}_{n} {\kern 1pt} \zeta \right|^{2} }  {\kern 1pt} {\kern 1pt} \left|d\zeta \right|=\frac{1}{2\pi i} \int \limits _{T}\frac{1-\left|a_{n} \right|^{2} }{(1-{\kern 1pt} {\kern 1pt} \overline{a}_{n} {\kern 1pt} \zeta )(\zeta -{\kern 1pt} {\kern 1pt} a_{n} {\kern 1pt} )}  {\kern 1pt} {\kern 1pt} d\zeta {\kern 1pt} {\kern 1pt} {\kern 1pt} =1,{\kern 1pt}  {\kern 1pt} {\kern
     1pt}$$
     $$\left|B_{n} (a_{n} )\right|\ge \delta (a)>0,$$

then

 $$\left\| h\right\| _{H^{1} } \le \left\| x\right\| _{l^{1} } /\delta (a){\kern 1pt} {\kern 1pt} {\kern 1pt} <\infty  {\kern 1pt} {\kern 1pt}{\kern 1pt} {\kern 1pt}  and  {\kern 1pt} {\kern 1pt} {\kern 1pt} {\kern 1pt} h\in H^{1} .$$

\

     Let     $f(z)=\sum \limits _{k\ge 0}\frac{\hat{h}(k)}{k+1}  {\kern 1pt} {\kern 1pt} z^{k+1} {\kern 1pt} {\kern 1pt} {\kern 1pt} {\kern 1pt} {\kern 1pt} {\kern 1pt} {\kern 1pt} {\kern 1pt} {\kern 1pt} {\kern 1pt} {\kern 1pt} {\kern 1pt} {\kern 1pt} {\kern 1pt} {\kern 1pt} {\kern 1pt} {\kern 1pt} {\kern 1pt} {\kern 1pt} ({\kern 1pt} {\kern 1pt} f'(z)=h(z){\kern 1pt} {\kern 1pt} {\kern 1pt} {\kern 1pt} )$ .

\

     Since   $f'\in H^{1} $  then  $f\in m_{\alpha } $  for all  $\alpha >0$  [3]   and

     $$f'(a_{k} )(1-\left|a_{k} \right|^{2} )=h(a_{k} )(1-\left|a_{k} \right|^{2} ){\kern 1pt} {\kern 1pt} {\kern 1pt} ={\kern 1pt} {\kern 1pt} x_{k} ,{\kern 1pt} {\kern 1pt} {\kern 1pt} {\kern 1pt} {\kern 1pt} {\kern 1pt} {\kern 1pt} {\kern 1pt} k=1,2,\ldots .$$

\

\

\

\

\noindent
{\small Department of Mathematics\\
        Technical University\\
        25, Tsanko Dijstabanov,\\
        Plovdiv, Bulgaria\\
        e-mail: peyyyo@mail.bg}

\end{document}